\documentclass[final]{amsart}
\usepackage{graphicx}
\newtheorem{theorem}{Theorem}%[section]
\newtheorem{theorem&definition}{Theorem and Definition}[section]
\newtheorem{lemma}{Lemma}%[section]
\newtheorem{corollaire}{Corollary}%[section]
\newtheorem{definition}{Definition}[section]

\newtheorem{proposition}{Proposition}%[section]
\newtheorem{claim}{Claim}[section]

\newcommand{\rg}{{\rm rg}\kern 0.12em}
\newcommand{\co}{{\rm co}\kern 0.12em}
\newcommand{\cl}{{\rm cl}\kern 0.12em}
\newcommand{\bd}{{\rm bd}\kern 0.12em}
\newcommand{\card}{{\rm card}\kern 0.12em}
\newcommand{\epi}{{\rm epi}\kern 0.12em}
\newcommand{\supp}{{\rm supp}\kern 0.12em}
\newcommand{\proj}{{\rm proj}\kern 0.12em}
\newcommand{\interieur}{{\rm int}\kern 0.12em}
\newcommand{\interior}{{\rm int}\kern 0.12em}
\newcommand{\sgn}{{\rm sgn}\kern 0.12em}
\newcommand{\dom}{{\rm Dom}\kern 0.12em}
\newcommand{\reach}{{\rm reach}\kern 0.12em}
\newcommand{\diam}{{\rm diam}\kern 0.12em}
\newcommand{\lip}{{\rm Lip}\kern 0.12em}
\newcommand{\vect}{{\rm vect}\kern 0.12em}
\newcommand{\R}{\mbox{\rm {I\hspace{-0.8mm}R}}}
\newcommand{\N}{\mbox{\rm I\hspace{-0.85mm}I\hspace{-1,2mm}N}}
\newcommand{\Ri}{{\rm I\kern -0.12em R}} % pour indices(Maths) 
\newcommand{\Ni}{{\rm I\kern -0.12em N}} % pour indices(Maths)

\def\Box{\ \vrule width0.2cm depth0cm height0.2cm}
\renewcommand{\a}[1]{\makebox[10mm][l]{\hspace{2mm} #1}}

\begin{document}
\setcounter{section} 0  \setcounter{assertion} {0}
%\baselineskip=7mm
%\abovedisplayskip=1mm
%\belowdisplayskip=3mm
%\abovedisplayshortskip=3mm
%\belowdisplayshortskip=3mm

\def\picture #1 by #2
(#3){\dimen0 = \hsize \advance\dimen0 by -#1 \divide \dimen0 by 2 \hskip \dimen0
\vbox to #2{\hrule width #1 height 0pt depth 0pt \vfill \special{picture #3}}}

\title{Volume of tubes, non polynomial behavior}
\author{Jean-Francois CROUZET$^{^{1}}$ and Marc-Olivier CZARNECKI$^{^{1}}$}
\maketitle
\setcounter{footnote}{1} 
\footnotetext{Institut de Mathematiques et Modelisation de Montpellier
  UMR 5149 CNRS, \'equipe ACSIOM (formerly Laboratoire d'Analyse Convexe), Universit\'e de Montpellier 2, place Eug\`ene Bataillon,
34095 Montpellier cedex 5, France. E-mail: marco@math.univ-montp2.fr,
crouzet@math.univ-montp2.fr

\thanks{We are indebted to Jacques Lafontaine, who gave us much
  information through his manuscript from the s\'eminaire
  Bourbaki~\cite{Lafontaine}, and to Luigi Ambrosio and Giovanni
  Belletini for advises and information, especially for telling us of the work of M. Z\"ahle.}
}
\begin{abstract}

The behavior of the volume of the tube $B(M,r)$, around a given
compact subset $M$ of $\R^n$, depending on $r$, is an old and
important question with relations to many fields, like differential
geometry, geometric measure theory, integral geometry, and also
probability and statistics. Federer~(1959) introduces the class of
sets with positive reach, for which the volume is given by a polynom in
$r$. For applications, in numerical analysis and statistics for
example, an ``almost'' polynomial behavior is of equal interest. We
exhibit an example showing how far to a polynom can be the volume of
the tube, for the simplest extension of the class of sets with
positive reach, namely the class of (locally finite) union of sets with
positive reach -satisfying a tangency condition- introduced by Z\"ahle
(1984).\end{abstract}

\section{Introduction}

Let $M$ be a compact subset of $\R^{n}$, and $r$ a non negative real number. Consider the tube (or $r$-neighborhood)
$$
\overline{B}(M,r)=\{x\in \R^{n}|\; d(M,x) \leq r \}
$$
an its volume
$$
{\mathcal{L}}_n(\overline{B}(M, r)),
$$ 
where ${\mathcal{L}}_n$ denotes the  Lebesgue measure on $\R^{n}$.\\

The volume ${\mathcal{L}}_n(\overline{B}(M, r))$ is polynomial in various useful cases. When the set $M$ is convex, and the corresponding polynomial is
called Steiner formula, named after the seminal work of Jakob Steiner~\cite[1840]{Steiner}.  When the set $M$ is a submanifold of class
$C^2$ (and $r$ smaller than a given $r_0$), it is given by Weyl's formula, named after  the paper of Hermann Weyl~\cite[1939]{Weyl}. Hebert Federer~\cite[1959]{Federer_Curvature_measures} introduced the sets of positive reach, in order to unify both approaches, and it is the widest known class to this day, for which the volume  ${\mathcal{L}}_n(\overline{B}(M, r))$ is a polynomial -for $ r$ small enough.\\

The motivation for finding a polynomial formula for  ${\mathcal{L}}_n(\overline{B}(M, r))$ first comes from the statistics -putting apart the early work of  Jakob Steiner. The seminal work of Herbert Hotelling ~\cite[1939]{Hotelling} describes and solves a class of statistical problems by giving a polynomial formula for (the tube around) curves in $\R^2$. It apparently motivated the celebrated generalization of H. Weyl, which is published right next to H. Hotelling's paper.\\

Let us mention the overwhelming interest of the (polynomial) volume of
tubes in probability and statistics. It  allows for large deviation
estimates, approximation of the tail probabilities, simultaneous
confidence and prediction bounds, construction of significance tests,....
We refer to the papers of Knowles and Siegmund~\cite{knowles_siegmund}, Johansen and
Johnstone~\cite{johansen_johnstone}, Naiman~\cite{naiman}, Sun~\cite{sun}... not being exhaustive of course.
Also Donnelly~\cite{Donnelly} related the volume of tubes to the PDEs -the heat equation.\\

The volume of the tube ${\mathcal{L}}_n(\overline{B}(M, r))$ is related to the $p$-dimensional Haussdorff measure of the set $M$ by the Minkowski content\footnote{Distinguish upper- and lower Minkowski content in general.} (on the right hand side)
$$
{\mathcal{H}}_p(M)=\lim_{r\to 0, r>0}{{\mathcal{L}}_n(\overline{B}(M, r))\over  \alpha(n-p)\> r^{n-p}},
$$
where $\alpha(i)={\mathcal{L}}_{i} (B_{\Ri^{i}}(0, 1))$, whenever $M$ is
$p$-rectifiable (see for example~\cite{Federer}).\\

In the various problems involving the volume
${\mathcal{L}}_n(\overline{B}(M, r))$, the fact that it is polynomial
is mainly used to obtain a rate of convergence. Our first motivation
was to numerically compute the perimeter ${\mathcal{H}}_{n-1}(\bd M)$ of
an n-dimensional set $M$, and it is quite obvious that a polynomial formula provides a polynomial rate of convergence (when $r\to 0$).\\

So, we only need an ``almost'' polynomial formula for the volume of
the tube  ${\mathcal{L}}_n(\overline{B}(M, r))$, say something like
$$
{\mathcal{L}}_n(\overline{B}(M, r))=P(r)+O(r^\lambda),
$$ 
in order to obtain the various estimates usually provided by an exact
polynomial formula -$P(r)$ being a polynomial, and the value of
$\lambda$ depending of the problem -for example, with $\lambda =n-p+1$
if the set $M$ is $p$-rectifiable, or  $\lambda =2$ if its boundary
$\bd M$ is $n-1$ rectifiable.\\

D.  Hug~\cite{Hug} and  J. Rataj~\cite{Rataj} show that, for a wide
class of compact subsets $M$ of $\R^{n}$ (locally finite union of sets with
positive reach in~\cite{Rataj}), having ${\mathcal{H}}_{n-1}$-almost everywhere one unit vector,\footnotemark
$$
{\mathcal{L}}_n(\overline{B}(M, r))={\mathcal{L}}_n(M)+
r{\mathcal{H}}_{n-1}(\bd M)+o(r).
$$ 

Let ${\mathcal{M}}$ be a class of compact n-dimensional subsets of
$\R^n$. For the purpose of efficiently (thus with a rate of
convergence) compute the perimeter of a set $M$ in  ${\mathcal{M}}$,
we would need the following result:

\subsection*{Conjecture (${\mathcal{M}}$)} {\it For every set $M\in
  {\mathcal{M}}$,
$$
{\mathcal{L}}_2(\overline{B}(M, r))={\mathcal{L}}_2(M)+
r{\mathcal{H}}_{n-1}(\bd M)+O(r^2).
$$
}
The theory of Federer shows the validity of Conjecture ($PR_n$), where
$PR_n$ is the set of compact sets of positive
reach, having ${\mathcal{H}}_{n-1}$-almost everywhere one unit vector.\footnotemark[3]\\
       
\footnotetext{We only make this restriction to obtain exactly
  ${\mathcal{H}}_{n-1}(\bd M)$ in the expansion, for a smooth reading.}

It is easy to give simple counterexamples for union of convex sets,
without any ``contact'' condition: consider the union of two tangent disks in $\R^2$
$$
M=\overline{B}((1,0),1)\cup \overline{B}((-1,0),1)
$$
Then\footnote{\begin{eqnarray*}
{\mathcal{L}}_2(\overline{B}(M, r))&=&2 (1+r)^2 \left(\pi-\arccos \left({1\over 1+r}\right)\right)+2(2r+r^2)^{1\over 2}\\
&=&2\pi +4\pi r-({8\sqrt{2}}/{3})r^{3/2}+O(r^ 2)\\
&=&
{\mathcal{L}}_2(M)+
r{\mathcal{H}}_{1}(\bd M)-({8\sqrt{2}}/{3})r^{3/2}+O(r^ 2).
\end{eqnarray*}
}
\begin{eqnarray*}
{\mathcal{L}}_2(\overline{B}(M, r))
&=&
{\mathcal{L}}_2(M)+
r{\mathcal{H}}_{1}(\bd M)+\alpha r^{3/2}+O(r^ 2).
\end{eqnarray*}

Extensions of Federer's Theory have been made in various
directions. In the Riemannian setting, by the work of Fu~\cite{Fu}. In
the Euclidean setting,  Z\"ahle~\cite{Zahle_Curvature_measures} considered the finite unions
of sets of positive reach satisfying a tangential condition -and gave
a polynomial formula for the ``volume'' of the tube. So did Cheeger,
M\"uller and Sch\"ader~\cite{CMS} in the case of piecewise linear sets
and R. Schneider~\cite{Schneider} for unions of convex sets.\\

But in both cases the corresponding ``value'' is not the volume
${\mathcal{L}}_n(\overline{B}(M, r))$. It is a modified volume, taking
into account the multiplicity of the normal vectors to the set
$M$. How big is the difference between the volume
${\mathcal{L}}_n(\overline{B}(M, r))$ and the modified value. Is it
small enough, for example of order $2$, to verify the above conjecture?\\

This is true in dimension 2. The purpose of this paper is to provide a
counterexample when the dimension is higher, for
which the volume ${\mathcal{L}}_n(\overline{B}(M, r))$ is ``far'' from
being a polynom, and which belongs the class $U_{PR}$ introduced by
Z\"ahle~\cite{Zahle_Curvature_measures}. This shows in
particular that Conjecture
($U_{{PR}_n}$) does not holds, for $n\geq 3$ (taking $U_{{PR}_n}$ to be the set of
compact  sets in $U_{PR}$), having ${\mathcal{H}}_{n-1}$-almost everywhere one unit vector\\

Our  counterexample (Theorem~\ref{contre-ex}) is the union of two
convex sets $M$ and $M'$ in $\R^3$, for which holds the nondegeneracy tangential
condition\footnote{Rataj and Z\"ahle later developed their theory
  without use of the nondegeneracy tangential
condition. However, as we point out thereafter, without this
condition and for our problem, the counterexample is obvious.}
 -defining the class $U_{PR}$ in~\cite{Zahle_Curvature_measures}:
$$
T(M\cap M',x)=T(M,x)\cap T(M',x),
$$
for every $x\in M\cap M'$. Since these are convex sets, the tangent cones are the usual ones, and this removes any hope of replacing the  tangent cones (namely Bouligand tangent cones) involved in the definition of the  class $U_{PR}$, by another tangent cone -in the quest of a quasi-polynomial formula.\\

%(or also two circles).\\

%Our counterexample is sharp, compared to the general behavior of the volume of the union of convex sets in $\R^3$ (Theorem~\ref{o(r)}).

\addtocounter{footnote}{1} \footnotetext{We let
$\R_{+}=\{x\in \R| x\geq 0\}$ and $\sgn \>x =x /|x|$ if $x\in \R\setminus \{0\}$. If
$x=(x_{1},...,x_{n})$ and $y=(y_{1},...,y_{n})$  belong to $\R^{n}$, we
 denote   $(x|y)= \sum_{i=1}^{n} x_{i} y_{i},$  the scalar product of $\R^{n}$,
$\|x\|=\sqrt{(x|x)},$ the Euclidean norm; we  denote $B(x,r)=\{y\in \R^{n}|\; \| x-y
\| < r \}$, $\overline{B}(x,r)=\{y\in \R^{n}|\; \| x-y \| \leq r \}$ and
$S(x,r)=\{y\in \R^{n}|\; \| x-y \| = r\}$. If $X\subset \R^{n},$ $Y\subset \R^{n}$,
and $x\in\R^{n}$, we let $d_X(x)=\inf_{y\in X} \|x-y\|$, $X \setminus Y =
\{x \in X | x \notin Y \}$ the set-difference of the sets $X$ and $Y$, $X+Y=\{x+y|x\in
X, y\in Y\}$, the sum of the sets $X $ and $Y$, $B(X,r)=X+B(0,r)$,
$\overline{B}(X,r)=X+\overline{B}(0,r)$, cl$X$, the closure of $X$,
$\interieur X$, the interior of $X$, $\bd X = $cl$X \setminus \interieur X$, the
boundary of $X$, $ \co X$, the convex hull of $X$.}

\section{Statement of the results}

\subsection{Main result}
Our main result shows that 
\begin{theorem}\label{contre-ex} For every integer  $N\geq 2$, there
  exists two convex compact sets $M$ and $M'$ in $\R^3$, such
  that\footnote{
In other words, the set $M\cup M'$ belongs to the class   $U_{PR}$ introduced by M. Z\"ahle~\cite{Zahle_Curvature_measures}.
Moreover,
$$
M\cup M'\subset \co\{ (0,0,0) ,(1,0,-1), (1,1,0), (1,-1,0)\}
$$
$$
 \co\{ (0,0,0) ,(1,0,-1), (1,1,0)\}\cup \co\{ (0,0,0) ,(1,0,-1),
 (1,-1,0)\}\subset M\cup M'
$$
}
\begin{eqnarray*} 
\forall x\in M\cap M',&& T(M\cap M',x)=T(M,x)\cap T(M',x),
\end{eqnarray*}
and the volume of  the tube $\overline{B}(M\cup M',r)$ satisfies%\footnote{PERSO FAIRE SAUTER LA CONSTANTE  $\frac{2}{3}$, PAR UNE HOMOTETHIE PAR EXEMPLE?}
$$
- \frac{16}{3} r^{1+\frac{1}{N}}+O(r^2) 
\leq {\mathcal{L}}_3(\overline{B}(M\cup M', r))- {\mathcal{L}}_3(M\cup M')- r{\mathcal{H}}_2(\bd(M\cup M'))\leq
- \frac{2}{3N}  r^{1+\frac{1}{N}}+O( r^{2})
.$$
\end{theorem}

Theorem~\ref{contre-ex} is proved in Section~\ref{contre-ex_proof}.
We can provide better asymptotic bounds (Theorem~\ref{contre-ex_1}
below) but the result is sharp in the sense that the``non polynomial''
part has to be negligible in front of $r$, as stated in the next
result, a special case of~\cite[Theorem 3.3]{Hug} and~\cite[Theorem 3]{Rataj}.

\begin{theorem} \label{o(r)} Let $ M$ and $M'$ in $\R^3$ be two convex compact subsets of $\R^3$, with nonempty interiors.%\footnote{If one does not assume that the sets have nonempty interiors, a corresponding result still holds by replacing.....} 
Then 
$$
{\mathcal{L}}_3(\overline{B}(M\cup M', r))- {\mathcal{L}}_3(M\cup M')- r{\mathcal{H}}_2(\bd(M\cup M'))\in o\left(r\right)
$$
\end{theorem}

Theorem~\ref{o(r)} is proved in Section~\ref{o(r)_proof}, for the sake
of completeness and with a simple self contained proof.

\begin{theorem}\label{contre-ex_1} 

For every function $\varepsilon: [0,1]\to \R_+$ of class $C^2$ such that
\begin{multline*}
\varepsilon(0)=0,\quad \lim_{r\rightarrow
  0}\frac{\varepsilon(r)}{\sqrt{r}}=+\infty,\quad 
\int_0^1\frac{\varepsilon(t)}{t}<+\infty,\\ \forall r\in (0,1],\quad
r^2\varepsilon''(r)-r\varepsilon'(r)+\varepsilon(r)>0
\end{multline*}
and the function $t\mapsto \frac{\varepsilon(t)}{t}$ is strictly convex and  decreasing,
there exists
    $\lambda\in\R$ and two convex compact sets $M$ and $M'$ in $\R^3$
such that:%\footnote{PERSO FAIRE SAUTER LA CONSTANTE  $2$, PAR UNE
	  %HOMOTETHIE PAR EXEMPLE? OU PAR LE CHOIX DE I(R)}
\begin{multline*}
-2r\varepsilon(r)+O(r^2) \geq {\mathcal{L}}_3(\overline{B}(M\cup M', r))- {\mathcal{L}}_3(M\cup
M')-r{\mathcal{H}}_2(\bd(M\cup M'))\\
\geq -16 r\int_0^r\frac{\varepsilon(t)}{t}{\rm d}t+\lambda r\varepsilon(r)+O(r^2).
\end{multline*}
\end{theorem}
%\footnote{PERSO I(R)
%\begin{theorem}
%For any strictly convex and strictly decreasing function $I$ of class
%$C^2$ on $(0,1]$ verifying  $I\in L^1([0,1])$, $\lim_{r\rightarrow
%0}rI(r)=0$, $\lim_{r\rightarrow 0}r^{1/2}I(r)=+\infty$ and
%$-\frac{I''(r)}{I'(r)}>\frac{1}{r},\forall r\in (0,1]$ there exists $\lambda\in\R$ and two convex compact sets $M$ and $M'$ in $\R^3$ such that:
%\begin{multline*}
%-2r^2I(r)+O(r^2) \geq {\mathcal{L}}_3(B(M\cup M', r))-
%{\mathcal{L}}_3(M\cup M')-r{\mathcal{H}}_2(\bd(M\cup M'))\\
%\geq -16 r\int_0^rI(t){\rm d}t+\lambda r^2 I(r)+O(r^2).
%\end{multline*}
%\end{theorem}
%}

Theorem~\ref{contre-ex_1} is proved in Section~\ref{contre-ex_1_proof}.\\

{\bf Remark.} As an example take  $\varepsilon>0$ and $I(r)=\frac{1}{r|\ln(r)|^{1+\varepsilon}}$ for
$r$ small enough.
This yields the estimate
\begin{multline*}
-2\frac{r}{|\ln(r)|^{1+\varepsilon}} +O(r^2) \geq {\mathcal{L}}_3(\overline{B}(M\cup M',
r))- {\mathcal{L}}_3(M\cup M')-r{\mathcal{H}}_2(\bd(M\cup M'))\\
\geq -\frac{16 r}{\varepsilon|\ln(r)|^{\varepsilon}}+\frac{\lambda
r}{|\ln(r)|^{1+\varepsilon}}+O(r^2).
\end{multline*}

{\bf Remark.} In view of the Steiner formula and of the proof of
Theorem~\ref{contre-ex}, one can replace the (Landau) functions
$O(r^2)$ by polynomials $a r^2+b r^3$ in the statement of
Theorem~\ref{contre-ex}, for $r$ small enough. By contrast, the result
in Theorem~\ref{contre-ex_1} is only asymptotic, see the end of Section~\ref{contre-ex_1_proof}.

\subsection{Relation to the class  $U_{PR}$}

In this section, we precisely recall the definitions of sets with
positive reach and $U_{PR}$ sets. \\

Theorem~\ref{contre-ex} is
a counterexample to a possible extension of a Steiner-Weyl type formula
-with ``small'' error, in the sense of Conjecture ($U_{{PR}_n}$)- for
the class of $U_{PR}$ sets, introduced
by~\cite{Zahle_Curvature_measures}. Sets in $U_{PR}$ are defined as
union of  sets with positive reach (introduced
by~\cite{Federer_Curvature_measures}),  satisfying a (nondegeneracy)
tangency condition.\\

Let $M\subset \R^n$ be nonempty. For $x\in \R^n$, the projection set of $x$ on $M$ is defined by:
$$ 
\proj_M(x)=\{ y\in M |   d(x,M)=\Vert y-x \Vert \} .
$$
The reach of $M$ is defined by:
$$\reach (M) =\sup\{ r>0 |\forall y \in B(M,r),\> \proj_M(y) \mbox{ reduces to a singleton}\}.$$
We let $\reach (\emptyset)=+\infty$.\\

We now recall the definition of sets with positive reach.

\begin{definition} [Federer~\cite{Federer_Curvature_measures}] A closed set  $M\subset \R^n$  is said to be of positive reach if $ \reach (M)>0$.
\end{definition}

{\bf Remark.} A closed set $M$ is of positive reach if it satisfies one of the following conditions (see~\cite{Federer_Curvature_measures}):\\
(i) $M$ is convex;\\
(ii) $M$ is a compact $C^2$ submanifold of $\R^n$, with or without a boundary.\\

In order to  generalize the Steiner-Weyl formula to the sets of
positive reach, Federer~\cite{Federer_Curvature_measures} builds a
general theory of curvature measures. The curvature measures give indeed the coefficient of the  Steiner-Weyl polynom.

\begin{theorem}
  [Federer,~\cite{Federer_Curvature_measures}]\label{polynom} Let
  $M\subset \R^n$ be compact\footnote{If one does not assume that $M$
    is bounded, a similar formula holds for
    ${\mathcal{L}}_n(\overline{B}(M, r)\cap \proj_M^{-1}(Q))$ where
    $Q$ is a bounded Borel subset of $\R^n$.} of positive reach.
 Then there are $(c_0,\dots ,c_n)\in \R  ^{n+1}$  such that, for every $r\in [0, \reach(M)]$:
$$
{\mathcal{L}}_n(\overline{B}(M, r))=\sum_{i=0}^{n} c_i r^i.
$$
\end{theorem}

Let $M\subset \R^n$ be nonempty, and $x\in M$. Then the Bouligand tangent cone to $M$ at $x$, denoted $T^B(M,x)$ is defined by:
$$
{T}^B(M,x)= \{v\in \R^{n}| \exists (\lambda _{k})_{k\in \Ni}, \lambda _{k}>0, \exists (y _{k})_{k\in \Ni}, y_{k}\in M, y_{k}\rightarrow x,
v=\lim_{k\rightarrow\infty} \lambda _{k}(y_{k}-x)\}.
$$

We now recall the definition of the class $U_{PR}$
in~\cite{Zahle_Curvature_measures}.

\begin{definition} \label{U_PR.def} [Z\"ahle~\cite{Zahle_Curvature_measures}] A closed set $M\subset \R^n$ is said to be $U_{PR}$ if there is a sequence $(M_k)_{k\in \Ni}$ of closed sets  with positive reach such that:\\
\a{(a)} $M=\cup_{k\in \Ni}M_k$;\\
\a{(b)} the sequence $(M_k)$ is locally finite, precisely, for every
$r>0$,\\
\a{}  the set $\{k\in \N | M_k\cap B(0,r)\ne \emptyset\}$ is finite;\\
\a{(c)} for every finite subset $I\subset \N$, $reach(\cap_{i\in I}M_{i})>0$\\
\a{(d)} for every finite subset $I\subset \N$, for every $x\in
\cap_{i\in I}M_{i}$,\\
\a{} $T^B(\cap_{i\in I}M_{i},x)=\cap_{i\in I} T^B(M_{i},x).$\\
We call $(M_k)_{k\in \Ni}$ a (not necessarily unique) decomposition of $M$.%\footnote{verif decomp minimale}
\end{definition}

{\bf Remark.} Note that, if $M$ is compact, the sequence $(M_k)_{k\in \Ni}$ clearly reduces to a finite family.\\

{\bf Remark.}  It is easy to notice that the (exact) Steiner-Weyl formula does
not hold in general, even in the  class $U_{PR}$ (with no need of the
counterexample in Theorem~\ref{contre-ex}!). In $\R^2$, consider $M=\cup_{i=1}^4 M_i$, where $M_i=\{(x,y)\in \R^2| (x-x_i)^2+(y-y_i)^2=1, |x|\geq 1, |y|\geq 1\}$ and $(x_1, y_1)=(1,0)$, $(x_2, y_2)=(0,1)$ $(x_3, y_3)=(-1,0)$, $(x_4, y_4)=(0,-1)$. Then $M$ belongs to the class $U_{PR}$.\footnote{A straightforward computations gives:
$$
{\mathcal{L}}_2(\overline{B}(M, r))=\left\{ \begin{array}{lll} 
\pi (r^2+8r)+a(r) & \mbox{ if }& r\in [0,1]\\
2\pi (r+1)^2+4(r^2-1)^{1/2}+2r^2(\arccos({(r^2-1)^{1/2}\over r})-\arccos({1\over r}))+a(r) & \mbox{ if }& r\in [1,\sqrt{2}]\\
2\pi(r+1)^2+a(r)+4& \mbox{ if }& r\in [\sqrt{2},+\infty)\\
\end{array}
\right. ,$$
$$\begin{array}{llll} \mbox{where } &a(r)&=&-2((r+1)^2-b(r))+\sqrt{2}(1+b(r))((r+1)^2-b(r))^{1/2}-4(r+1)^2\arccos({1+b(r)\over 2(r+1)}),\\ \mbox{and }& b(r)&=&(2(r+1)^2-1)^{1/2}\end{array}$$
\footnote{$a(r)=-2((r+1)^2-(2(r+1)^2-1)^{1/2})+2\sqrt{2}(1+(2(r+1)^2-1)^{1/2})((r+1)^2-(2(r+1)^2-1)^{1/2})^{1/2}-8(r^2+1)\arccos({1+(2(r+1)^2-1)^{1/2}\over 2(r+1)})$.}
}\\

M. Z\"ahle studies~\cite{Zahle_Curvature_measures} the Steiner Weyl formula for $U_{PR}$ sets, by defining a modified volume of sets which essentially takes into account the multiplicity of normal cones, like Cheeger, M\"uller, Schrader
\cite{CMS} for  piecewise linear spaces, and R. Schneider~\cite{Schneider} for unions
of convex sets. To make it short, she obtains
a polynomial formula for this modified volume  of the tube
$\overline{B}(M,r)$, by adding the volumes of the tubes
$\overline{B}(M_k,r)$ of the decomposition.\\

 It is obvious
(above remark) that the (exact) polynomial formula will not hold in general
for the true volume ${\mathcal{L}}_n(\overline{B}(M, r))$  of the tube, even for
$U_{PR}$ sets. But one may wonder, in the light of the later results
of Hug~\cite{Hug} and Rataj~\cite{Rataj}, how far from a polynomial will the
volume behave.\\

In dimension 2, one obtains a Steiner Weyl type formula, with an extra
term $O(r^2)$. In dimension greater than $3$, 
   Theorem~\ref{contre-ex} shows the possibly bad
behavior of the volume of the tube, thus hindering any hope to extend
the (true, without multiplicity) Steiner Weyl formula -with ``small'' error- for $U_{PR}$ sets. Since our counterexample is the union of two convex sets, it also shows that replacing the Bouligand tangent cone by another tangent cone in the definition of $U_{PR}$ sets would not make any difference, since all tangent cones coincide in the convex case.

\section{Proof of Theorems~\ref{contre-ex} and~\ref{contre-ex_1}, a counterexample in
  dimension 3 \label{contre-ex_proof}}

We first prove a more general but more technical result,
which also helps to understand the values of the bounds. We deduce
Theorem~\ref{contre-ex} in Section~\ref{contre-ex_proof_1}, and Theorem~\ref{contre-ex_1}  in Section~\ref{contre-ex_1_proof}.

\begin{theorem} \label{contre-ex_2} Let $\varphi:[0,1]\to \R_+$ and
  $\psi :[0,1]\to \R_+$ be two strictly convex functions of class
  $C^2$ such that $\psi(0)=\varphi(0)=0$, $\psi(1)=\varphi(1)=1$, $\psi'(0)=0$,
  $\lim_{t\rightarrow 0}\frac{\varphi(t)}{\psi(t)}=0$ and
  $\frac{\varphi(t)}{\varphi'(t)}<\frac{\psi(t)}{\psi'(t)}, \forall t\in(0,1]$. 

We define 
$$
M =\co\Big( \left\{(t,0, -\psi(t))|t\in [0,1]\} \cup \{(t,\varphi(t), 0)|t\in [0,1]\right\}\Big),
$$
and let $M' $ be the symmetric of $M $ with respect to the plane $\{(x_1,x_2,x_3)| x_2=0\}$, precisely, $M' =\co\Big( \{(t,0, -\psi(t))|t\in [0,1]\} \cup \{(t,-\varphi(t), 0)|t\in [0,1]\}\Big)$.
The sets $M$ and $M'$ are compact and convex,
$$
\forall x\in M\cap M',\quad T(M\cap M',x)=T(M,x)\cap T(M',x)
$$
and
\begin{multline*}
-2r^2\int_{\varphi^{-1}(r)}^{1}\frac{\psi(t)}{\varphi(t)}{\rm d}t +O(r^2)\\
\geq {\mathcal{L}}_3(\overline{B}(M\cup M', r))- {\mathcal{L}}_3(M\cup M')-r{\mathcal{H}}_2(\bd(M\cup M'))\\
\geq -\left( 2 \alpha^{1/2} r^2
\int_{\rho(r)}^{1}\frac{\psi'(t)}{\varphi'(\theta(t))}{\rm d}t
+2 r \int_0^{\rho(r)} \psi(t)dt\right)+O(r^2),
\end{multline*}
where $\rho(r)\in [0,1]$, and  $\theta(t)$ is the unique solution of
$$
\theta(t)-\frac{\varphi(\theta(t))}{\varphi' (\theta(t))}=t-\frac{\psi(t)}{\psi'(t)}
$$
hence the function $\theta$ is continuous, and where
$$
\alpha=\max_{t\in [0,1]}\quad 1+
\frac{\varphi'(\theta(t))^2}{\psi'(t)^2}\left(1+ \psi(t)^2\right)
$$
\end{theorem}

\begin{corollaire} \label{x^p,x^q} For every integers $p$ and $q$ such
  that $q>p+1\geq 3$, the functions $\varphi(t)=t^q$ and $\psi(t)=t^p$ satisfy the assumptions of
  Theorem~\ref{contre-ex_2}, and yield, for the corresponding sets $M$
  and $M'$
$$
-\frac{2}{q-p-1} r^{1+\frac{p+1}{q}}+O(r^2)
 \geq {\mathcal{L}}_3(\overline{B}(M\cup M', r))- {\mathcal{L}}_3(M\cup M')-
r{\mathcal{H}}_2(\bd(M\cup M'))\geq  - c(p,q)  r^{1+\frac{p+1}{q}}+O(r^2)
$$
where
$$
c(p,q)=2\left(
\frac{p^2}{q^2}\left(\frac{q-1}{q}\frac{p}{p-1}\right)^{2q-2}+2
\right)^\frac{p+1}{2q}\>\left(\frac{1}{q-p-1}+\frac{1}{p+1}\right)
$$

\end{corollaire}

\begin{figure}[htb]
\centering
\includegraphics[width=0.5\textwidth]{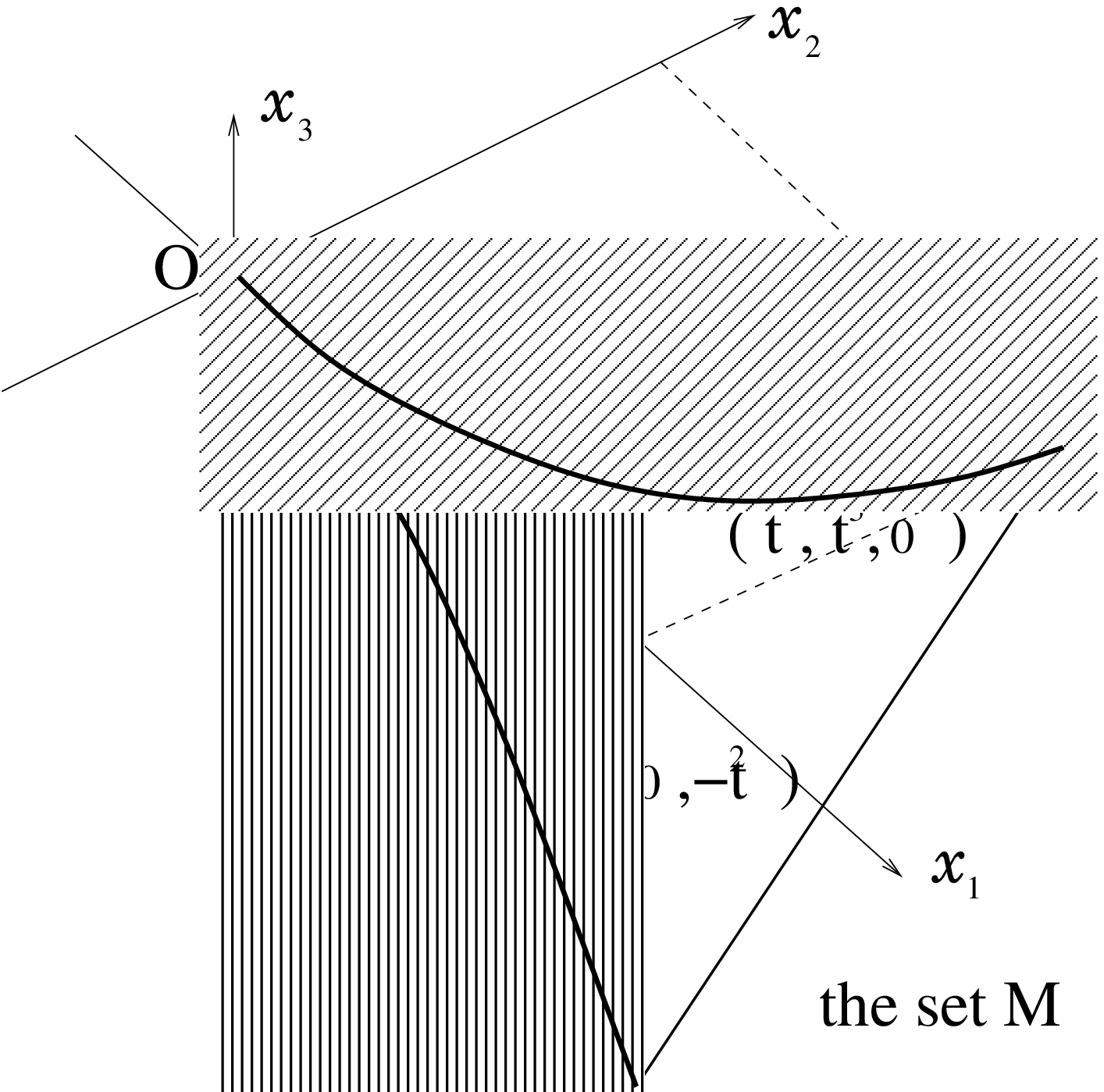}
\caption{The set $M$ with $\psi(t)=t^2$ and $\varphi(t)=t^3$}\label{f}
\end{figure}

\subsection{Proof of Theorem~\ref{contre-ex_2}} 

The proof of Theorem~\ref{contre-ex_2} goes in four steps. First we check the tangency condition, i.e.,  the set $M \cup M' $ belongs to the class $U_{PR}$.
Afterward we introduce the set which provokes the deviation from the polynomial term and we specify its volume. Then we give an upper bound of the volume ${\mathcal{L}}_3(\overline{B}(M \cup M' , r))$, and finally we give a lower bound of this volume.

\subsubsection{The set $M \cup M' $ belongs to the class $U_{PR}$} 
The sets $M $ and $M' $ are clearly compact and convex. Note that 
$$
M = \co (C \cup C_1),
$$
where 
\begin{eqnarray*}
C&=& \co  \{(t,0,-\psi(t))|t\in [0,1]\}\\
C_1&=&\co \{(t,\varphi(t), 0)|t\in [0,1]\}
\end{eqnarray*}
 Then $M
\cap M' =C$. One easily checks the tangency condition 
$$
\forall x\in M\cap M', T(M\cap M',x)=T(M,x)\cap T(M',x)
$$
by using  $\psi'(0)=0$ for the calculus at the origin, hence the set $M \cup M' $ belongs to the class $U_{PR}$.$\Box$

\subsubsection{The set $A(r)$}

Define 
$$
A(r)=\overline{B}(M , r)\cap \overline{B}(M' , r)\setminus\overline{B}(M \cap M' , r).
$$

We specify the volume of the set $A(r)$, and thus reduce the proof of  Theorem~\ref{contre-ex_2} mainly to the estimation of ${\mathcal{L}}_3 (A(r))$.

\begin{lemma}\label{vol_A(r)} Let $M$ and $M'$ be two compact convex subsets of $\R^3$, with nonempty interiors. Then
$$
{\mathcal{L}}_3 (A(r))= -{\mathcal{L}}_3(\overline{B}(M\cup M', r))+{\mathcal{L}}_3(M\cup M')+ r{\mathcal{H}}_2(\bd(M\cup M'))+ O(r^2).
$$
\end{lemma}

{\bf Proof of Lemma~\ref{vol_A(r)}.} From the following partition of $\overline{B}(M \cup M', r)$
\begin{multline*}
\overline{B}(M \cup M' , r)=\overline{B}(M , r)\setminus \left(\overline{B}(M , r)\cap \overline{B}(M' , r)\right)
\,\, \sqcup \,\, \overline{B}(M' , r)\setminus \left(\overline{B}(M , r)\cap \overline{B}(M' , r)\right)\\
\sqcup\>  \left(\overline{B}(M , r)\cap \overline{B}(M' , r)\right)
\end{multline*}
and from the partition
$$
\overline{B}(M , r)\cap \overline{B}(M' , r)= A(r) \sqcup \overline{B}(M \cap M' , r)
$$
we deduce
$$
{\mathcal{L}}_3(\overline{B}(M\cup M', r))= {\mathcal{L}}_3(\overline{B}(M, r))+{\mathcal{L}}_3(\overline{B}( M', r))-{\mathcal{L}}_3(\overline{B}(M\cap M', r))-{\mathcal{L}}_3 (A(r))
$$
From Steiner's formula,
\begin{eqnarray*}
 {\mathcal{L}}_3(\overline{B}(M, r))&=&{\mathcal{L}}_3(M)+r{\mathcal{H}}_2(\bd M)+ O(r^2)\\
{\mathcal{L}}_3(\overline{B}(M', r))&=&{\mathcal{L}}_3(M')+r{\mathcal{H}}_2(\bd M')+ O(r^2)
\end{eqnarray*}
For the set $M\cap M'$, depending on its dimension, Steiner's formula gives%\footnote{PERSO A VERIFIER, REFERENCER, ET A UNIFIER EVENT.}
\begin{eqnarray*}
{\mathcal{L}}_3(\overline{B}(M\cap M', r))&=& {\mathcal{L}}_3(M\cap M')+r{\mathcal{H}}_2(\bd (M\cap M'))+ O(r^2) \hbox{ if } \dim (M\cap M')=3\\
{\mathcal{L}}_3(\overline{B}(M\cap M', r))&=& 2 r{\mathcal{H}}_2(\bd (M\cap M'))+ O(r^2) \hbox{ if } \dim (M\cap M')=2\\
{\mathcal{L}}_3(\overline{B}(M\cap M', r))&=& O(r^2) \hbox{ if } \dim (M\cap M')\leq 1
\end{eqnarray*}
From the partition 
$$M \cap M'= M\setminus (M\cap M') \sqcup  M'\setminus (M\cap M') \sqcup  M\cap M',$$
$${\mathcal{L}}_3(M\cap M')={\mathcal{L}}_3(M)+{\mathcal{L}}_3(M')- {\mathcal{L}}_3(M\cap M').$$
For the boundaries, we use the following partitions
\begin{eqnarray*}
\bd(M\cup M')&=&\bd M\setminus M' \sqcup \bd M'\setminus M\sqcup ( \bd M \cap \bd M')\setminus \interior (M\cup M')\\
\bd M&=&\bd M\setminus M' \sqcup ( \bd M \cap \bd M')\setminus \interior (M\cup M')\sqcup (\bd M\cap\interior (M\cup M')) \\
\bd M'&=&\bd M'\setminus M \sqcup ( \bd M \cap \bd M')\setminus \interior (M\cup M')\sqcup (\bd M'\cap\interior (M\cup M'))
\end{eqnarray*}
and the decomposition
$$
\bd(M\cap M')= ( \bd M \cap \bd M')\setminus \interior (M\cup M') \cup  (\bd M\cap\interior (M\cup M')) \cup  (\bd M'\cap\interior (M\cup M')) 
$$
We now consider the different possibilities for $\dim (M\cap M')$.

\paragraph{Case $\dim (M\cap M')=3$}
A Hahn-Banach separation argument yields $\bd M\cap \bd M'\cap\interior (M\cup M')=\emptyset$. Hence the above decomposition of $\bd(M\cup M')$ is a partition hence
$${\mathcal{H}}_2(\bd (M\cup M'))={\mathcal{H}}_2(\bd M)+{\mathcal{H}}_2(\bd M')-{\mathcal{H}}_2(\bd (M\cap M')).$$
Combined with the above Steiner's formulas, this proves  Lemma~\ref{vol_A(r)}.

\paragraph{Case $\dim (M\cap M')\leq 2$} When $\dim (M\cap M')\leq 2$, and since $M$ and $M'$ are convex with nonempty interiors,
\begin{eqnarray*}
\bd M\cap\interior (M\cup M') &= & \bd M'\cap\interior (M\cup M')\\
M\cap M'=\bd(M\cap M')&= &\bd M\cap \bd M'
\end{eqnarray*}
Moreover ${\mathcal{H}}_2( ( \bd M \cap \bd M')\setminus \interior
(M\cup M'))=0$ (distinguish the cases $\dim (M\cap M')\leq 1$ and
$\dim (M\cap M')=2$, and use the fact that $\dim M=\dim M'=3$). Hence
$${\mathcal{H}}_2(\bd (M\cup M'))={\mathcal{H}}_2(\bd M)+{\mathcal{H}}_2(\bd M')- 2 {\mathcal{H}}_2(\bd (M\cap M')).$$
Combined with the above Steiner's formulas, this proves  Lemma~\ref{vol_A(r)}.$\Box$

\subsubsection{Upper bound of the volume ${\mathcal{L}}_3(\overline{B}(M \cup M' , r))$}

In view of Lemma~\ref{vol_A(r)}, it is sufficient to give a lower bound of the volume of $A(r)$. 

\begin{lemma}\label{minoration} For $r>0$, let 
\begin{multline*}
\check{A}(r)=\left\{(x_1,x_2,x_3)\in (0,1]\times \R^2| x_2\leq 0, x_1\ge\varphi^{-1}(r),\phantom{\frac{\psi^2(x_1)}{\varphi^2(x_1)}} 
\right.\\
\left. -\frac{\varphi(x_1)}{\psi(x_1)}x_2-\psi(x_1)\leq x_3\leq x_2\frac{\psi(x_1)}{\varphi(x_1)}-\psi(x_1)+r\sqrt{1+\frac{\psi^2(x_1)}{\varphi^2(x_1)}} \right\}
\end{multline*}
Then
\begin{eqnarray}
\label{minoration_1}&\check{A}(r)\subset \overline{B}(M ,r)\cap \overline{B}(M' ,r) \cap \{(x_1,x_2,x_3)|x_2\leq 0\};&\\
\label{minoration_2}&\check{A}(r)\setminus \overline{B}(\{(t,0, -\psi(t))|t\in [0,1]\},r)
\subset  A(r) \cap \{(x_1,x_2,x_3)|x_2\leq 0\}.&
\end{eqnarray}
\end{lemma} 

Let us first admit the lemma. Since ${\mathcal{L}}_3(A(r))=2
{\mathcal{L}}_3(A(r)\cap \{(x_1,x_2,x_3)|x_2\leq 0\})$, then
${\mathcal{L}}_3(A(r))\geq 2 {\mathcal{L}}_3(\check{A}(r))-
{\mathcal{L}}_3(\overline{B}(\{(t,0, -\psi(t))|t\in
[0,1]\},r))$. Noting that 
$${\mathcal{L}}_3(\check{A}(r))=\int_{\varphi^{-1}(r)}^1r^2\frac{\psi(t)}{\varphi(t)}{\rm d}t
$$ 
and that, for $r$ small enough,
$$
{\mathcal{L}}_3(\overline{B}(\{(t,0, -\psi(t))|t\in [0,1]\},r)=c
r^2+d r^3, 
$$ where $c={\mathcal{H}}_1(\{(t,0,-\psi(t))|t\in [0,1]\})$ and $d={4 \over 3}\pi$, we obtain
$$
{\mathcal{L}}_3(A(r))\geq 2r^2\int_{\varphi^{-1}(r)}^{1}\frac{\psi(t)}{\varphi(t)}{\rm d}t -c r^2-d r^3.\Box
$$

{\bf Proof of Lemma~\ref{minoration}.} {\it Proof of \rm (\ref{minoration_1}).} Consider an element $x=(x_1,x_2,x_3)$ in the set $\check{A}(r)$. 
The proof consists in checking that its distance to $M $ is less or
equal to $r$ by considering its projection $p$ on the line containing
$(x_1,0,-\psi(x_1))$ and $(x_1, \varphi(x_1), 0)$. Let%\footnote{PERSO $r>0$ donc $x_1>0$ ...}    
$$
p=\left|\begin{array}{lll} p_1&=&x_1\vspace{1mm}\\ 
p_2&=&\displaystyle  x_2-\frac{1}{1+\frac{\varphi^2(x_1)}{\psi^2(x_1)}}\left(x_2-\frac{\varphi(x_1)}{\psi(x_1)}x_3-\varphi(x_1)\right)\vspace{1mm}\\
p_3&=&\displaystyle x_3+\frac{\frac{\varphi(x_1)}{\psi(x_1)}}{1+\frac{\varphi^2(x_1)}{\psi^2(x_1)}}\left(x_2-\frac{\varphi(x_1)}{\psi(x_1)}x_3-\varphi(x_1)\right)\vspace{1mm}\end{array}\right. 
$$

We now check that $p\in \co \{(x_1,0,-\psi(x_1)), (x_1, \varphi(x_1),
0)\}\subset M $. Indeed, 
$$
p=\theta (x_1,0,-\psi(x_1))+(1-\theta)(x_1, \varphi(x_1), 0),$$ with 
$$
\theta=\frac{1}{1+\frac{\varphi^2(x_1)}{\psi^2(x_1)}}\left(\frac{\varphi^2(x_1)}{\psi^2(x_1)}-x_2\frac{\varphi(x_1)}{\psi^2(x_1)}-\frac{x_3}{\psi(x_1)}\right).
$$
By assumption, $-\frac{\varphi(x_1)}{\psi(x_1)}x_2-\psi(x_1)\leq
x_3$, hence $\theta\leq 1$. By assumption, $r\leq
\varphi(x_1)$\footnote{since $\varphi$ is increasing}, $x_2\leq 0$, and $ x_3\leq x_2\frac{\psi(x_1)}{\varphi(x_1)}-\psi(x_1)+r\sqrt{1+\frac{\psi^2(x_1)}{\varphi^2(x_1)}}$, hence
\begin{eqnarray*}
 \theta\left( {1+\frac{\varphi^2(x_1)}{\psi^2(x_1)}}\right) &\geq &\frac{\varphi^2(x_1)}{\psi^2(x_1)}-x_2\frac{\varphi(x_1)}{\psi^2(x_1)}-x_2\frac{1}{\varphi(x_1)}+1-\varphi(x_1)\sqrt{\frac{1}{\psi^2(x_1)}+\frac{1}{\varphi^2(x_1)}}\\
 &\geq & 1+\frac{\varphi^2(x_1)}{\psi^2(x_1)}- \sqrt{1+\frac{\varphi^2(x_1)}{\psi^2(x_1)}}\geq 0
\end{eqnarray*}
Note that
$$
\|p-x\|=\frac{\vert x_2
  -\frac{\varphi(x_1)}{\psi(x_1)}x_3-\varphi(x_1)\vert}{\sqrt{1+\frac{\varphi^2(x_1)}{\psi^2(x_1)}}}.
$$
By assumption, noticing that $\varphi$ is increasing, $x_3\leq
  x_2\frac{\psi(x_1)}{\varphi(x_1)}-\psi(x_1)+r\sqrt{1+\frac{\psi^2(x_1)}{\varphi^2(x_1)}}$, hence $ \frac{x_2 -\frac{\varphi(x_1)}{\psi(x_1)}x_3-\varphi(x_1)}{\sqrt{1+\frac{\varphi^2(x_1)}{\psi^2(x_1)}}}\geq -r$. By assumption, $x_2\leq 0$, and  $-\frac{\varphi(x_1)}{\psi(x_1)}x_2-\psi(x_1)\leq x_3 $. Hence
 $ x_2 -\frac{\varphi(x_1)}{\psi(x_1)}x_3-\varphi(x_1) \leq  -\frac{\varphi(x_1)}{\psi(x_1)}x_3-\varphi(x_1)\leq \frac{\varphi^2(x_1)}{\psi^2(x_1)}x_2\leq 0$. Finally, 
$$d(M ,x)\leq\|p-x\|\leq r.$$
Now take  
$$
p'=\left|\begin{array}{c} p_1\\-p_2\\p_3\end{array}\right. .
$$
Since $p\in \co \{(x_1,0,-\psi(x_1)), (x_1, \varphi(x_1), 0)\}$, then $p'$ belong to  the symmetrical set $\co \{(x_1,0,-\psi(x_1)), (x_1, -\varphi(x_1), 0)\}$ hence $p'\in M'$. Since $x_2\leq 0$ and $p_2 \geq 0$ (since $p\in M $), we obtain  that $(x_2-p'_2)^2=(x_2+p_2)^2\leq(x_2-p_2)^2$, hence that $d(M' ,x)\leq \|p'-x\|\leq \|p-x\|\leq r$.$\Box$\\

{\it Proof of \rm (\ref{minoration_2}).} Consider an element $x=(x_1,x_2,x_3)$ in the
left-hand side set. We now prove that, if $x\notin A(r)$, then $x\in
\overline{B}(\{(t,0, -\psi(t))|t\in [0,1]\},r)$. Indeed, since $x\in
(\overline{B}(M ,r)\cap \overline{B}(M' ,r))\setminus A(r)$ by~(\ref{minoration_1}), from the
definition of the set $A(r)$, $x\in \overline{B}(M \cap M'
,r)=\overline{B}(C,r)$. Let $\proj_{C}(x)=(q_1,0,q_3)$. If
$\proj_{C}(x)\notin \{(t,0, -\psi(t))|t\in [0,1]\}$, then
$x-\proj_{C}(x)=(-\lambda,\mu,-\lambda)$ for some $(\lambda, \mu)\in
\R_+\times \R$, and $\lambda =0$ if $q_3\ne-q_1$. Noticing that  $q_3+\psi(q_1)\leq 0$, the case  $\lambda
=0$ implies that $q_1=x_1$ and $q_3=x_3$ hence
$x_3+\psi(x_1)\leq 0$. Hence, the definition of $\check{A}(r)$ yields
$x_2=0$ and $x_3+\psi(x_1)\geq 0$, hence $q_3+\psi(q_1)\geq 0$. As  $q_3+\psi(q_1)\leq 0$, we deduce $q_3=\psi(q_1)=0$, hence $q_1=0$ (by the strict convexity of $\psi$) and $\proj_{C}(x)=(0,0,0)$, a
contradiction. The case $q_3=-q_1$ implies that $x_1+x_3=-2\lambda\leq
0$, a contradiction with the definition of the left-hand side set,
which implies that $-\frac{\varphi(x_1)}{\psi(x_1)}x_2-\psi(x_1)\leq
x_3$, hence $x_1+x_3\geq \psi(x_1)+x_3\geq
-\frac{\varphi(x_1)}{\psi(x_1)}x_2\geq 0$ (since $\psi$ is
convex with $\psi(0)=0$, $\psi(1)=1$, hence $x_1\geq\psi(x_1)$, and since we work in the half space $\{(x_1,x_2,x_3)|x_2\leq 0\}$).$\Box$

\subsubsection{Lower bound of the volume ${\mathcal{L}}_3(\overline{B}(M \cup M' , r))$}
We first verify the definition of the function $\theta$ and the real number $\alpha$.

\begin{claim}\label{theta_alpha} There exists a  continuous function  $\theta: [0,1]\to [0,1]$ such that, for every $t$,  $\theta(t)$ is the unique solution of
$$
\theta(t)-\frac{\varphi(\theta(t))}{\varphi' (\theta(t))}=t-\frac{\psi(t)}{\psi'(t)}
.$$
Besides, we can define
$$
\alpha=\max_{t\in [0,1]}\quad 1+
\frac{\varphi'(\theta(t))^2}{\psi'(t)^2}\left(1+ \psi(t)^2\right)
$$
\end{claim}

{\bf Proof of Claim~\ref{theta_alpha}.}
Note that if $f$ is (as $\psi$ and $\varphi$) a strictly convex
function of class $C^2$ on $[0,1]$ that satisfies $f(0)=0$, we have
$0<\frac{f(t)}{f'(t)}<t$ on $(0,1]$ and $t\mapsto
t-\frac{f(t)}{f'(t)}$ is continuous and strictly increasing on
$[0,1]$. Hence the reciprocal function of $t\mapsto
t-\frac{\varphi(t)}{\varphi' (t)}$ is defined on
$\left[0,1-\frac{\varphi(1)}{\varphi' (1)}\right]$ and is continuous.
Similarly,  $t\mapsto t-\frac{\psi(t)}{\psi'(t)}$ is a continuous
strictly increasing function mapping $[0,1]$ on
$\left[0,1-\frac{\psi(1)}{\psi'(1)}\right]\subset
\left[0,1-\frac{\varphi(1)}{\varphi'(1)}\right]$, as by assumption
$\frac{\varphi(1)}{\varphi'(1)}<\frac{\psi(1)}{\psi'(1)}$.
Hence the solution $\theta(t)$ of the equation
$\theta(t)-\frac{\varphi(\theta(t))}{\varphi'(\theta(t))}=t-\frac{\psi(t)}{\psi'(t)}$
is unique and continuous. \\

We have moreover by assumption
$-\frac{\varphi(t)}{\varphi'(t)}+\frac{\psi(t)}{\psi'(t)}>0$, and we
recall that  $-t+\frac{\psi(t)}{\psi'(t)}<0$. As $\theta(t)$ is the
unique point annulating the continuous function
$x\mapsto x-\frac{\varphi(x)}{\varphi'(x)}-t+\frac{\psi(t)}{\psi'(t)}$, we
  have $0<\theta(t)<t$ by the mean value theorem.
Thus, as $\psi'$ and $\varphi'$ are clearly non negative functions,
and $\varphi'$ is increasing, we have $0\leq\lim_{t\rightarrow
  0}\frac{\varphi'(\theta(t))}{\psi'(t)}\leq \lim_{t\rightarrow
  0}\frac{\varphi'(t)}{\psi'(t)}=0$, given that $\lim_{t\rightarrow
  0}\frac{\varphi'(t)}{\psi'(t)}=\lim_{t\rightarrow
  0}\frac{\varphi(t)}{\psi(t)}=0$ by l'H\^opital's rule. Hence
$t\mapsto 1+\frac{\varphi'(\theta(t))^2}{\psi'(t)^2}\left(1+\psi(t)^2\right)$
is continuous on $[0,1]$ and attains its supremum $\alpha$.$\Box$\\

Back to the lower bound of the volume ${\mathcal{L}}_3(\overline{B}(M
\cup M' , r))$,  it is now sufficient to give an upper bound of the
volume of $A(r)$ in view of Lemma~\ref{vol_A(r)},

\begin{lemma}\label{majoration}
Let 
$$
\widehat{A}(r)=\left\{(x_1,x_2,x_3)| x_2\leq 0,-\psi(x_1)\leq x_3,\> -x_2+ \frac{\varphi'(\theta(x_1))}{\psi'(x_1)}(x_3+\psi(x_1))\leq \alpha^{1/2} r\right\}.
$$
Then 
$$
A(r)\cap \{(x_1,x_2,x_3)|x_1\in [0,1],\> x_2\leq 0\}\subset \widehat{A}(r).
$$
\end{lemma}

Admitting the lemma, we obtain 
$${\mathcal{L}}_3(A(r)\cap [\rho(r),1]\times \R^2)\leq
2{\mathcal{L}}_3(\widehat{A}(r)\cap [\rho(r),1]\times
\R^2)=2\alpha^{1/2} r^2 \int_{\rho(r)}^{1}\frac{\psi'(t)}{\varphi'(\theta(t))}{\rm d}t.$$
In view of Lemma~\ref{majoration}, we notice that 
$$
A(r)\cap [0, \rho(r)]\times \R^2\subset \{(x_1,x_2,x_3)|
0\leq x_1\leq \rho(r), -r\leq x_2\leq r,-\psi(x_1)\leq x_3\leq
r\}.
$$
If $x_1\leq 0$ and $x_2\leq 0$ (if $x_2\geq 0$, consider $d(M',x)$) we deduce by a direct computation that
$$
d(M,x)\geq d(\co\{(0,0,0), (1,0,0),
(1,1,0), (1,0,-1)\},x)\geq d(\co\{(0,0,0), (1,0,-1)\},x)\geq d(M \cap
M' ,x) $$ which contradicts the definition of the set $A(r)$. Hence
$$
A(r)\cap (-\infty,0]\times \R^2=\emptyset
$$
Finally,
$$
A(r)\cap [1, +\infty)\times \R^2\subset [1, 1+r] \times [-r,r] \times  [-1-r,1].
$$
and we obtain
$$
{\mathcal{L}}_3(A(r))\leq  2 \alpha^{1/2} r^2
\int_{\rho(r)}^{1}\frac{\psi'(t)}{\varphi'(\theta(t))}{\rm d}t
+ 2 r \int_0^{\rho(r)} \psi(t)dt +2 r^2 \rho(r)+2 r^2(1+2r)   .\Box$$

{\bf Proof of Lemma~\ref{majoration}.} Consider an element $x=(x_1,x_2,x_3)$ in the set $A(r)\cap \{(x_1,x_2,x_3)|x_2\leq 0\}$. First note that $-r\leq x_2$. Indeed, since $M\subset\{(x_1,x_2,x_3)|x_2\geq 0\}$, we have $r\geq d(M,x)\geq d(\{(x_1,x_2,x_3)|x_2\geq 0\},x)=\|x-(x_1,0,x_3)\|=-x_2$.\\
 
We first prove that $-\psi(x_1)\leq x_3$. Assume that it is not true. If $-x_1\leq x_3$, then 
$$
(x_1,0, x_3)\in \co \{(x_1,0,-\psi(x_1)), (x_1,0,-x_1)\}\subset C=M \cap M' .
$$
 Hence $d(M \cap M' ,x)\leq |x_2|\leq r$, which contradicts the fact that $x\in A(r)$. Now assume that $x_3<-x_1$. Note that $n=(-1,0,-1)$ strictly separates $x$ and $M $. Indeed, $(n|x)=-x_1-x_3>0$, and, for every $\lambda \in [0,1]$, $(n|(\lambda ,0, -\psi(\lambda)) )=-\lambda +\psi(\lambda)\leq 0$, $(n|(\lambda ,\varphi(\lambda) ,0) )=-\lambda \leq 0$, hence $(n|\Lambda)\leq 0$ for every $\Lambda\in M $. In particular, $(n|\proj_{M }(x))\leq 0$, which implies that $(n|x-\proj_{M }(x))>0$ hence that $n\notin T(M ,\proj_{M }(x))$. Note also that $n'=(-1,-2,-1)$ strictly separates $x$ and $M $, $(n'|x)=-x_1-2x_2-x_3> 0$, hence that $n'\notin T(M ,\proj_{M }(x))$. We now show that this implies that $\proj_{M }(x)\in C=M \cap M' $, which contradicts the fact that $x\in A(r)$. 
Recall that
\begin{eqnarray*}
C&=& \co  \{(t,0, -\psi(t))|t\in [0,1]\}=\{(y_1,0,y_3)| -y_1\leq y_3\leq -\psi(y_1)\}\\
C_1&=&\co \{(t,\varphi(t), 0)\>|\>\, t\in [0,1]\}\>=\{(z_1,z_2,0)|\varphi(z_1)\leq z_2\leq z_1\}.
\end{eqnarray*}
Since $\proj_{M }(x)\in M =\co(C\cup C_1)$, there are $y\in C$, $z\in C_1$ and $\theta \in [0,1]$ such that $\proj_{M }(x)=\theta y + (1-\theta) z$. 
 If $\theta=1$ or $z= (0,0,0)$, then clearly $\proj_{M }(x)\in C $. Now assume that $\theta<1$ and $z\ne (0,0,0)$.  If $z\notin \co\{(0,0,0),(1,1,0)\}$, for $\varepsilon $ small enough, $z+\varepsilon n\in M $ (for example, $z+ {z_1-z_2\over 2} n= z_2 (1,1,0)+ {z_1-z_2\over 2} (1,0,-1) \in \co \{(0,0,0), (1,1,0), (1,0, -1)\}$), hence $\proj_{M }(x)+\varepsilon (1-\theta)n\in M $, which implies that $n\in T(M ,\proj_{M }(x))$, a contradiction. If $z\in \co\{(0,0,0),(1,1,0)\}\setminus \{(0,0,0)\}$, then $z+\varepsilon n'\in M $ (for example, $z+ {z_1\over 2} n'= {z_1\over 2} (1,0,-1) \in \co \{(0,0,0), (1,1,0), (1,0, -1)\}$), hence $\proj_{M }(x)+\varepsilon (1-\theta)n'\in M $, which implies that $n'\in T(M ,\proj_{M }(x))$, a contradiction.\\

We now prove that $-x_2+ \frac{\varphi'(\theta(x_1))}{\psi'(x_1)}(x_3+\psi(x_1))\leq \alpha^{1/2} r\ $. This a consequence of the following claim.

\begin{claim} \label{normale} $$n''=\left(\varphi'(\theta(x_1)) , -1,
  \frac{\varphi'(\theta(x_1))}{\psi'(x_1)}\right)\in N\Big(M , (x_1, 0,-\psi(x_1))\Big).$$
\end{claim}

Admitting the claim, we deduce
\begin{eqnarray*}
(x-(x_1, 0, -\psi(x_1))\,|\,n'')&=&(x-\proj_{M }(x)+\proj_{M }(x)-(x_1, 0, -\psi(x_1))\,|\,n'')\\
&\leq &\|x-\proj_{M }(x)\| \, \|n''\| \leq r\|n''\|
\end{eqnarray*}
Hence 
$$ -x_2+\frac{\varphi'(\theta(x_1))}{\psi'(x_1)}(x_3+\psi(x_1))\leq
\sup_{t\in [0,1]}\quad
\left(1+\frac{\varphi'(\theta(t))^2}{\psi'(t)^2}\left(1+\psi(t)^2\right)\right)^{1/2} r = \alpha^{1/2} r.
\Box$$

{\bf Proof of Claim~\ref{normale}.} Consider $y=(y_1,0,y_3)\in
C$. Then, in view of the definition of $C$,
$$
(n''|y-(x_1, 0, -\psi(x_1)))=n_1^{''}(y_1-x_1)+n_3^{''}(y_3+\psi(x_1))\leq
n_1^{''}(y_1-x_1)+n_3^{''}(-\psi(y_1)+\psi(x_1))
$$
 Since $n_1^{''}=\psi'(x_1)n_3^{''}$,  and since
$\psi$ is convex,
$$(n''|y-(x_1, 0,-\psi(x_1))) \leq n_3^{''}(-\psi(y_1)+\psi(x_1)+(y_1-x_1)\psi'(x_1))\leq 0.$$
Now consider $z=(z_1,z_2,0)\in C_1$. From the
definition of $C_1$, $z_2\geq \varphi(z_1)$, hence $$(n''|z-(x_1,0,-\psi(x_1)))\leq
(z_1-x_1)\varphi'(\theta(x_1))-\varphi(z_1)+\psi(x_1)\frac{\varphi'(\theta(x_1))}{\psi'(x_1)}.$$   
The study of the concave function,  $\lambda \mapsto
(\lambda-x_1)\varphi'(\theta(x_1))-\varphi(\lambda)+\psi(x_1)\frac{\varphi'(\theta(x_1))}{\psi'(x_1)}$ shows that it
attains  its maximum on $\R_+$ at the point $\theta(x_1)$ where its value is zero. Hence  $(n''|z-(x_1, 0, -\psi(x_1)))\leq 0$. Since $M =\co(C\cup C_1)$, it ends the proof of the claim.$\Box$

\subsection{Proof of Corollary~\ref{x^p,x^q}} 

It is immediate to check that the functions $\varphi(t)=t^q$ and $\psi(t)=t^p$ satisfy the assumptions of
  Theorem~\ref{contre-ex_2}, and provide the explicit value
  $\theta(t)=\frac{q}{q-1}\frac{p-1}{p}t$. Hence, we get 
\begin{eqnarray*}
\int_{\rho(r)}^{1}\frac{\psi'(t)}{\varphi'(\theta(t))}{\rm d}t&=&\frac{p}{q}\left(\frac{q-1}{q}\frac{p}{p-1}\right)^{q-1}\frac{\rho(r)^{p+1-q}-1}{q-p-1}\\
\int_0^{\rho(r)} \psi(t)dt&=&\frac{\rho(r)^{p+1}}{p+1}\\
\alpha&=& 1+ 2\frac{q^2}{p^2}\left(\frac{q}{q-1}\frac{p-1}{p}\right)^{2q-2}.
\end{eqnarray*}
Take
$$
\rho(r)=\left(
\frac{p^2}{q^2}\left(\frac{q-1}{q}\frac{p}{p-1}\right)^{2q-2}+2
\right)^\frac{1}{2q}\> r^\frac{1}{q}
$$
which belongs to $[0,1]$ for $r$ small enough and
which yields the lower bound:
$$
{\mathcal{L}}_3(\overline{B}(M\cup M', r))- {\mathcal{L}}_3(M\cup
M')-r{\mathcal{H}}_2(\bd(M\cup M'))\geq -  c(p,q)r^{1+\frac{p+1}{q}}+O( r^{2})
$$
For the upper bound, 
$$
\int_{\varphi^{-1}(r)}^{1}\frac{\psi(t)}{\varphi(t)}{\rm d}t=\frac{r^{\frac{p+1}{q}-1}-1}{q-p-1} 
$$
leads to 
$$
-\frac{2}{q-p-1} r^{1+\frac{p+1}{q}}+O(r^2) \geq {\mathcal{L}}_3(\overline{B}(M\cup M', r))- {\mathcal{L}}_3(M\cup M')-r{\mathcal{H}}_2(\bd(M\cup M'))
$$ which ends the proof of the corollary.$\Box$

\subsection{Proof of Theorem~\ref{contre-ex} \label{contre-ex_proof_1}}

Take $p=2$ and $q=3 N$ in  Corollary~\ref{x^p,x^q}.$\Box$

\subsection{Proof of Theorem~\ref{contre-ex_1}\label{contre-ex_1_proof}}

Take $\psi(t)=t^2$, $\rho(r)=2\varphi^{-1}(r)$ and
$I(r)=\int_{\varphi^{-1}(r)}^{1}\frac{\psi(t)}{\varphi(t)}{\rm d}t$ in
Theorem~\ref{contre-ex_2}.
The definition of $\theta(t)$ and the convexity of $\psi$ imply
$\frac{\psi'(t)}{\varphi'(\theta(t))}\leq\frac{\psi}{\varphi}(\theta(t))$.
Since    $\theta(t)\geq t-\frac{\psi(t)}{\psi'(t)}$, and since
 $\frac{\psi}{\varphi}$ is decreasing, we get
$\int_{2\varphi^{-1}(r)}^{1}\frac{\psi'(t)}{\varphi'(\theta(t))}{\rm d}t\leq
\int_{2\varphi^{-1}(r)}^{1}\frac{\psi}{\varphi}\left(t-\frac{\psi(t)}{\psi'(t)}\right){\rm d}t$.
The choice of $\psi(t)=t^2$ yields
$$
\int_{2\varphi^{-1}(r)}^{1}\frac{\psi'(t)}{\varphi'(\theta(t))}{\rm d}t \leq
\int_{2\varphi^{-1}(r)}^{1}\frac{\psi}{\varphi}\left(\frac{t}{2}\right){\rm d}t \leq 2I(r).
$$
We have moreover  $\int_0^{2\varphi^{-1}(r)} \psi(t){\rm
d}t=8\int_0^{\varphi^{-1}(r)} \psi(t){\rm d}t$.
From the definition of $I(r)$ a simple derivation and integration by
parts gives $\int_0^{\varphi^{-1}(r)} \psi(t){\rm d}t=\int_0^r-tI'(t){\rm
d}t=\int_0^rI(t){\rm d}t-rI(r)$ (as by assumption $\lim_{r\rightarrow
0}rI(r)=0$).
The estimation proved in Theorem~\ref{contre-ex_2} reads now as expected
(with $\lambda=16-4\alpha$).
\begin{multline*}
-2r^2I(r)+O(r^2) \geq {\mathcal{L}}_3(\overline{B}(M\cup M', r))- {\mathcal{L}}_3(M\cup
M')-r{\mathcal{H}}_2(\bd(M\cup M'))\\
\geq -\left( 4 \alpha r^2 I(r)+16 r\left(\int_0^rI(t){\rm
d}t-rI(r)\right)\right) +O(r^2).
\end{multline*} 

We now have to check that $\varphi$ given by the relation
$I(r)=\int_{\varphi^{-1}(r)}^{1}\frac{\psi(t)}{\varphi(t)}{\rm d}t=\frac{\varepsilon(r)}{r}$
satisfies the assumptions of Theorem~\ref{contre-ex_2}.
The relation $\int_0^{\varphi^{-1}(r)} \psi(t){\rm d}t=\int_0^r\frac{\varepsilon(t)}{t}{\rm d}t-\varepsilon(r)$
and the choice of $\psi(t)=t^2$
yields
 $$
\varphi^{-1}(r)=3\left(\int_0^r\frac{\varepsilon(t)}{t}{\rm d}t-\varepsilon(r)\right)^{1/3}.
$$
Hence $\varphi(0)=0$. A simple derivation proves that
$\varphi$ is a strictly increasing function.
Moreover, the condition  $\lim_{r\rightarrow
  0}\frac{\varepsilon(r)}{\sqrt{r}}=+\infty$
and the convexity of $t\mapsto \frac{\varepsilon(t)}{t}$ prove that $\lim_{t\rightarrow
  0}\frac{\varphi(t)}{t^2}=0$ (by a Taylor expansion and using the
above expression of $\varphi^{-1}(r)$).
The condition $r^2\varepsilon''(r)-r\varepsilon'(r)+\varepsilon(r)>0$ is sufficient to
verify that $\frac{\varphi(t)}{\psi(t)}=\frac{\varphi(t)}{t^2}$ is
strictly increasing on $[0,1]$, that is
$\frac{\varphi(t)}{\varphi'(t)}<\frac{\psi(t)}{\psi'(t)},\forall t\in(0,1]$.
As by assumption $t\mapsto \frac{\varepsilon(t)}{t}$  is a strictly convex function of class $C^2$ on
$(0,1]$ and $\frac{\varphi}{\psi}$ is strictly increasing, a double derivation
of
$I(r)=\int_{\varphi^{-1}(r)}^{1}\frac{\psi(t)}{\varphi(t)}{\rm d}t$ proves that
$\varphi^{-1}$ is strictly concave, and hence that $\varphi$ is strictly
convex. 
Finally, the condition $\varphi(1)=1$ is obtained by modifying the
function $\varphi$ on the interval $[\beta, 1]$ for a sufficiently
small $\beta$, keeping the properties of $\varphi$ (strict convexity,
$\frac{\varphi}{\psi}$ strictly increasing)\footnote{An easy way is to
  take $\beta$, such that $\varphi(\beta)<1$
  replace $\varphi$ with $t\mapsto \frac{\varphi(\beta)-1}{\beta^2-1}
  t^2+\frac{\beta^2-\varphi(\beta)}{\beta^2-1}$ on  $[\beta, 1]$, and
 carefully smoothing the new function around $\beta$.}, noticing that the
asymptotic behavior when $r\to 0$ remains unchanged.$\Box$

\section{Proof of Theorem~\ref{o(r)} \label{o(r)_proof} }

In view of Lemma~\ref{vol_A(r)}, it is sufficient to prove the following lemma.

\begin{lemma}\label{vol_A(r)=o(r)} Let $M$ and $M'$ be two compact convex subsets of $\R^3$. Then
$$
{\mathcal{L}}_3 \left(A(r)\right)\in o(r)
$$
\end{lemma}

{\bf Proof of Lemma~\ref{vol_A(r)=o(r)}.}
From the following partition of the set $\bd\left(M\cup M'\right)$
$$
\bd\left(M\cup M'\right)=\bd M\setminus M' \sqcup \bd M'\setminus M\sqcup \left( \bd M \cap \bd M'\right)\setminus \interior \left(M\cup M'\right)
$$
we write, using the elementary fact $ \overline{B}\left(A\cup B,r\right)= \overline{B}\left(A,r\right)\cup  \overline{B}\left(B,r\right)$,
\begin{multline*}
\overline{B}\left(\bd\left(M\cup M'\right),r\right)=\overline{B}\left(\bd M\setminus M',r\right) \>\cup \>\overline{B}\left( \bd M'\setminus M,r\right)\\
\cup \>\overline{B}\left( \left( \bd M \cap \bd M'\right)\setminus \interior \left(M\cup M'\right),r\right)
\end{multline*}
Hence, using the elementary formula ${\mathcal{L}}_3\left(A\cup B\cup C\right)={\mathcal{L}}_3\left(A\cup C\right)+{\mathcal{L}}_3\left(B\cup C\right)-{\mathcal{L}}_3\left( C\right)-{\mathcal{L}}_3\left(A\cap B\setminus A\cap B\cap C\right)$
\begin{multline*}
{\mathcal{L}}_3\left(\overline{B}\left(M\cup M', r\right)\right)={\mathcal{L}}_3\left(\overline{B}\left(\bd M\setminus M' \quad   \cup \quad   \left( \bd M \cap \bd M'\right)\setminus \interior \left(M\cup M'\right) ,r\right)\right)\\
\hfill+{\mathcal{L}}_3\left(\overline{B}\left(\bd M'\setminus M \quad   \cup \quad   \left( \bd M \cap \bd M'\right)\setminus \interior \left(M\cup M'\right) ,r\right)\right)\\
\hfill- {\mathcal{L}}_3\left(\overline{B}\left( \left( \bd M \cap \bd M'\right)\setminus \interior \left(M\cup M'\right) ,r\right)\right)\\
-{\mathcal{L}}_3\left(\left(\overline{B}\left(\bd M\setminus M',r\right)\cap \overline{B}\left( \bd M'\setminus M,r\right)\right)\setminus \overline{B}\left( \left( \bd M \cap \bd M'\right)\setminus \interior \left(M\cup M'\right) ,r\right)\right)
\end{multline*}
Now notice that each set in the above decomposition of  the set $\bd\left(M\cup M'\right)$ is $2$-rectifiable, as part of the boundary of a convex subset of $\R^3$, or union of two such sets. Hence by~\cite[Theorem 3.2.39]{Federer} their Minkowski content equals the $2$-dimensional Hausdorff measure of their closure. Noticing additionally that the set $\bd M\setminus M'\sqcup \left( \bd M \cap \bd M'\right)\setminus \interior \left(M\cup M'\right)$ is closed, we deduce
\begin{multline*}
{\mathcal{H}}_2\left(\bd\left(M\cup M'\right)\right)=\lim_{r\to 0}\frac{1}{2r}{\mathcal{L}}_3\left(\overline{B}\left(M\cup M', r\right)\right)\hfill
\end{multline*}
\begin{multline*}
{\mathcal{H}}_2\left(\bd M\setminus M'\sqcup \left( \bd M \cap \bd M'\right)\setminus \interior \left(M\cup M'\right)\right)=\\ 
\lim_{r\to 0}\frac{1}{2r}{\mathcal{L}}_3\left(\overline{B}\left(\bd M\setminus M'\sqcup \left( \bd M \cap \bd M'\right)\setminus \interior \left(M\cup M'\right),r\right)\right)
\end{multline*}
\begin{multline*}
{\mathcal{H}}_2\left(\bd M'\setminus M\sqcup \left( \bd M \cap \bd M'\right)\setminus \interior \left(M\cup M'\right)\right)=\\
\lim_{r\to 0}\frac{1}{2r}{\mathcal{L}}_3\left(\overline{B}\left(\bd M'\setminus M\sqcup \left( \bd M \cap \bd M'\right)\setminus \interior \left(M\cup M'\right),r\right)\right)
\end{multline*}
\begin{multline*}
{\mathcal{H}}_2\left(\left( \bd M \cap \bd M'\right)\setminus \interior \left(M\cup M'\right)\right)=\\
\lim_{r\to 0}\frac{1}{2r}{\mathcal{L}}_3\left(\overline{B}\left( \left( \bd M \cap \bd M'\right)\setminus \interior \left(M\cup M'\right),r\right)\right)
\end{multline*}
Remark that
\begin{multline*}
{\mathcal{L}}_3\left(\bd\left(M\cup M'\right)\right)={\mathcal{L}}_3\left(\bd M\setminus M' \>   \cup \>   \left( \bd M \cap \bd M'\right)\setminus \interior \left(M\cup M'\right)\right)\\
\hfill +
{\mathcal{L}}_3\left(\bd M'\setminus M \>  \cup \>   \left( \bd M \cap \bd M'\right)\setminus \interior \left(M\cup M'\right)\right)\\
-
{\mathcal{L}}_3\left(\left( \bd M \cap \bd M'\right)\setminus \interior \left(M\cup M'\right) \right)
\end{multline*}
and that
$$
A\left(r\right) \subset \overline{B}\left(\bd M\setminus M',r\right)\cap \overline{B}\left( \bd M'\setminus M,r\right)\setminus \overline{B}\left( \left( \bd M \cap \bd M'\right)\setminus \interior \left(M\cup M'\right),r\right)
$$
and deduce
$$
\lim_{r\to 0}\frac{{\mathcal{L}}_3 \left(A\left(r\right)\right)}{2r}=0.\Box
$$
Let us just point out that we may not have ${\mathcal{H}}_2\left(\cl\left(\bd M\setminus M'\right)\right)={\mathcal{H}}_2\left(\bd M\setminus M'\right)$ without additional assumptions, which prevents us to use a more straightforward decomposition and a rough bound.

\end{document}